\documentclass[12pt]{amsart}
\usepackage{pdfsync}

\theoremstyle{plain}
\newtheorem{de}{Definition}[section]
\newtheorem{thm}{Theorem}[section]
\newtheorem{rem}[thm]{Remark}
\newtheorem{cor}[thm]{Corollary}
\newtheorem{prop}[thm]{Proposition}
\newtheorem{lem}[thm]{Lemma}

\begin{document}
\newcommand{\R}{{\mathbb R}}
\newcommand{\He}{{\mathbb H}}
\newcommand{\N}{{\mathbb N}}
\newcommand{\Rn}{{\R^n}}
\newcommand{\Rp}{{\R_+^n}}
\newcommand{\beq}{\begin{equation}}
\newcommand{\eeq}{\end{equation}}
\newcommand{\beqs}{\begin{equation*}}
\newcommand{\eeqs}{\end{equation*}}
\newcommand{\ieqa}{\begin{eqnarray}}
\newcommand{\eeqa}{\end{eqnarray}}
\newcommand{\ieqas}{\begin{eqnarray*}}
\newcommand{\eeqas}{\end{eqnarray*}}
\newcommand{\al}{\alpha}
\newcommand{\lam}{\lambda}
\newcommand{\lams}{\lambda^+_\al}
\newcommand{\lamso}{\lambda^-_\al}
\newcommand{\upl}{u^+_\al}
\newcommand{\umi}{u^-_\al}
\newcommand{\Om}{\Omega}
\newcommand{\Oms}{\overline{\Omega}}
\newcommand{\om}{\omega}
\newcommand{\p}{\partial}
\newcommand{\ep}{\epsilon}
\newcommand{\xs}{\overline{x}}
\newcommand{\ys}{\overline{y}}
\newcommand{\xa}{x_\al}
\newcommand{\gt}{\widetilde{g}}
\renewcommand{\dim}{\noindent\textbf{Proof.} }
\newcommand{\dims}{\noindent\textbf{Proof} }
\newcommand{\finedim}{{\unskip\nobreak\hfil\penalty50
   \hskip2em\hbox{}\nobreak\hfil\mbox{$\Box$ \qquad}
   \parfillskip=0pt \finalhyphendemerits=0\par\medskip}}

\def\A{{\mathcal A}}
\def\S{{\mathcal S}}
\def\H{{\mathcal H}}
\def\M{{\mathcal{M}_{a,A}^+}}
\def\T{{\mathcal T}}
\def\I{{\mathcal I}}
\def\F{{\mathcal F}}
\def\J{{\mathcal J}}
\def\E{{\mathcal E}}
\def\P{{\mathcal P}}
\def\HH{{\mathcal H}}
\def\V{{\mathcal V}}
\def\B{{\mathcal B}}

\font\teneufm=eufm10 \font\seveneufm=eufm7 \font\fiveeufm=eufm5
\newfam\eufmfam
\textfont\eufmfam=\teneufm \scriptfont\eufmfam=\seveneufm
\scriptscriptfont\eufmfam=\fiveeufm
\def\eufm#1{{\fam\eufmfam\relax#1}}
\def\Def{{\eufm D}}
\def\Man{{\eufm M}}

\newcommand{\average}{{\mathchoice {\kern1ex\vcenter{\hrule
height.4pt width 6pt depth0pt} \kern-11pt}
{\kern1ex\vcenter{\hrule height.4pt width 4.3pt depth0pt}
\kern-7pt} {} {} }}
\newcommand{\ave}{\average\int}
\title{A Neumann eigenvalue problem for fully nonlinear operators.}
\author{Isabeau Birindelli, Stefania Patrizi}
\date{}
\maketitle
\begin{center} {\em Dedicated to prof. Louis Nirenberg for his $85^{\rm{th}}$ birthday.}\end{center}
\section{Introduction}
In this introduction and in the rest of the paper we quote some
works of Louis Nirenberg that are used explicitly in order to give
the right definitions and to prove the results; but the influence
of his research, here and in all the papers both the authors have
written, goes well beyond the citations. His mathematical ideas
have been very important for us, specially for the first named
author, but his teaching of how to approach mathematical problems
has been as important. We are happy to have this opportunity to
thank him for his generosity.

In this paper, for $\Omega$ a $C^2$ bounded domain of $\R^n$ and
for any $\alpha>0$, we consider the eigenvalue problem:
\begin{equation}\label{000}
\begin{cases}
 \M(D^2u)+\lambda u=0 & \text{in} \quad\Om, \\
 \frac{\p u}{\p
\overrightarrow{n}}=\alpha u & \text{on} \quad\partial\Om, \\
 \end{cases}
\end{equation}
where $\M$ is the Pucci operator, i.e.
$\displaystyle\M(M)=\sup_{0<aI\leq\sigma\leq AI}\rm{tr}(\sigma M)
$.

It is useless to emphasize the importance of the concept of
eigenvalue for the understanding of the structural properties of
the solutions both for linear and non linear equations. The
pioneering work of Berestycki, Nirenberg and Varadhan \cite{BNV}
has open the way to enlarge this fundamental concept to non linear
operators. Indeed, even if they treat linear equations, their
theory is very well adapted to fully nonlinear operators and
viscosity solutions being based primarily on the use of the
maximum principle. This has been done by many and in many
different contests, let us mention the works of Armstrong,  Busca,
Demengel, Juutinen, Ishii, Quaas, Sirakov, Yoshimura and the
authors of this note (\cite{SA,bd,BEQ,IY,J,p2,QS1}). It should be
mentioned that P.-L. Lions in \cite{PLL}, with a completely
different approach, first introduces what he called
demi-eigenvalues. Indeed when the operator is not odd with respect
to the Hessian (as is the case of the Pucci operators),
eigenvalues corresponding to positive eigenfunctions or to
negative eigenfunctions may not coincide and one could interpret
these two eigenvalues as a "splitting" of the eigenvalue.

The eigenvalue problem for Robin boundary conditions associated
with a fully-nonlinear operator  was already treated in \cite{p2}.
The novelty here is that we consider $\alpha>0$ which is the
"wrong sign" in the sense that the boundary conditions are not
"proper", see e.g. \cite{cil}. The boundary source  and the
reaction-diffusion equation are somehow in competition.

In analogy to \cite{BNV} we define the eigenvalues in the
following way:

\beqs\begin{split}\lams:=\sup\{&\lam\in \R\;|\;\exists\,v>0 \text{
on $\Oms$ bounded viscosity supersolution of }\\& \M(D^2v)+\lam
v=0 \text{ in } \Om,\, \frac{\p v}{\p \overrightarrow{n}} =\al
v\text{ on }\partial \Om \},
\end{split}\eeqs

\beqs\begin{split}\lamso:=\sup\{&\lam\in \R\;|\;\exists\,v<0
\text{ on $\Oms$ bounded viscosity subsolution of }\\&
\M(D^2v)+\lam v=0\text{ in } \Om,\, \frac{\p v}{\p
\overrightarrow{n}} =\al v\text{ on }\partial \Om \}.
\end{split}\eeqs

The first step is to prove that there exists $\upl>0$ and $\umi<0$
solutions of (\ref{000}) when respectively $\lambda=\lams$ and
$\lambda=\lamso$ (Proposition \ref{esistautof}).  We shall also
prove that below these eigenvalues there are solutions of the
equation with a forcing term $f(x)$ as long as the $f$ has the
right sign, i.e. $f\leq 0$ below $\lams$ and $f\geq0$ below
$\lamso$.

We are mainly interested in the asymptotic behavior with respect
to $\alpha$ of the eigenvalues. When  $\alpha\rightarrow 0$,
$\lams$ and $\lamso$ tend to $0$ which is the principal eigenvalue
of the pure  Neumann boundary problem
\begin{equation*}
\begin{cases}
 \M(D^2u)+\lambda u =0& \text{in} \quad\Om, \\
 \frac{\p u}{\p
\overrightarrow{n}}=0 & \text{on} \quad\partial\Om, \\
 \end{cases}
\end{equation*}

But our main goal is to study the behavior when
$\alpha\rightarrow+\infty$, this is done in our main

\begin{thm}\label{mainthm} The following limits hold:
\beq\label{lamsbehavior}\lim_{\al\rightarrow+\infty}\frac{\lams}{-\al^2}=A,\eeq
\beq\label{lamsobehavior}\lim_{\al\rightarrow+\infty}\frac{\lamso}{-\al^2}=a.\eeq
\end{thm}
Interestingly this asymptotic behavior emphasizes the "splitting"
of the eigenvalue. In the linear case, i.e. when $a=A=1$ and the
Pucci operator is nothing else but the Laplacian, this problem was
treated  in \cite{lz} by Lou and Zhu with a variational approach.
Very recently Daners and Kennedy \cite{dk} have proved that this asymptotic behavior is valid for the whole spectrum.

We also prove that for any $K\subset\subset\Omega$, the normalized
eigenfunctions $\upl$ and $\umi$ satisfy
$$\|\upl\|_{L^\infty(K)}\rightarrow 0\quad\mbox{ and }\quad \|\umi\|_{L^\infty(K)}\rightarrow 0\quad\text{as }\al\rightarrow+\infty.$$
So that the eigenfunctions tend to concentrate on the point of the
boundary where they reach the sup or the inf.

 The idea of the proof of Theorem \ref{mainthm} which somehow follows the line adopted in \cite{lz}, is the following: first we establish that $\upl$ reaches its maximum on the boundary and then we perform a blow up around this point.

 Then a key tool  will be a Liouville result in the half space (Theorem \ref{halfspacethm}).  Precisely we prove that for $\gamma>A$ (respectively $\gamma>a$) there are no bounded subsolutions (respectively supersolutions) of
$$ \left\{%
\begin{array}{ll}
    \M(D^2 u)- \gamma u=0 & \hbox{in }\Rp, \\
    -\frac{\p u}{\p x_n}=u & \hbox{on } \p\Rn.\\
\end{array}%
\right.
$$
 that are positive  (repectively negative) somewhere.
In \cite{lz} the analogous result for the Laplacian is proved
using the construction of sub and super solutions in the flavor of
what is done in \cite{bcn}.
 Let us mention here that it would be interesting to extend the results of Berestycki, Caffarelli, Nirenberg \cite{bcn} in half spaces,
  to this class of fully-nonlinear operators and to these boundary conditions.

 Lipschitz estimates
up to the boundary will be required in the proofs of both the
existence results and the asymptotic behavior. These estimates
which are interesting in their own right, are established here
using an argument inspired by \cite{il} (see also Barles and Da
Lio \cite{bdl} and Milakis and Silvestre \cite{ms} ).

In the whole paper the fully-nonlinear operator considered is the
Pucci operator $\M$, but, mutatis mutandis, parallel results can
be stated for the Pucci operator ${\mathcal{M}_{a,A}^-}$ defined
by $\displaystyle{\mathcal{M}_{a,A}^-}(M)=\inf_{0<aI\leq\sigma\leq
AI}\rm{tr}(\sigma M) $.

\section{Preliminary results}
Let us recall the definition of viscosity sub and supersolution of
the Neumann problem associated to a general elliptic operator
$F:\overline{\Om}\times\R\times\R^n \times
\emph{S(n)}\rightarrow\R$. Here $\emph{S(n)}$ is the space of
symmetric matrices on $\R^n$, equipped with the usual ordering. We
denote by $USC(\Oms)$ (resp., $LSC(\Oms)$) the set of upper
(resp., lower) semicontinuous functions on $\Oms$. Let
$f:\Oms\rightarrow\R$, $g:\p\Om\times\R\rightarrow\R$.
\begin{de} A function $u\in USC(\Oms)$ (resp., $u\in LSC(\Oms)$ )  is
called \emph{viscosity subsolution} (resp., \emph{supersolution})
of
\begin{equation*}
\begin{cases}
 F(x,u,Du,D^2u)=f(x) & \text{in} \quad\Om, \\
 \frac{\p u}{\p
\overrightarrow{n}}=g(x,u) & \text{on} \quad\partial\Om, \\
 \end{cases}
 \end{equation*}if the following conditions hold
\begin{itemize}
\item[(i)] For every $x_0\in \Om$, for any $\varphi\in C^2(\overline{\Om})$, such that $u-\varphi$ has a local maximum (resp., minimum)
at $x_0$ then
$$F(x_0,u(x_0),D\varphi(x_0),D^2\varphi(x_0))\geq \,(\text{resp., } \leq\,)\,f(x_0).$$
\item[(ii)]For every $x_0\in \partial\Om$, for any $\varphi\in C^2(\overline{\Om})$, such that $u-\varphi$ has a local maximum (resp., minimum)
at $x_0$ then
$$-(F(x_0,u(x_0),D\varphi(x_0),D^2\varphi(x_0))-f(x_0))\wedge
\left(\frac{\p \varphi}{\p
\overrightarrow{n}}(x_0)-g(x_0,u(x_0))\right)\leq 0$$(resp.,
$$-(F(x_0,u(x_0),D\varphi(x_0),D^2\varphi(x_0))-f(x_0))\vee
\left(\frac{\p \varphi}{\p
\overrightarrow{n}}(x_0)-g(x_0,u(x_0))\right)\geq
0).$$\end{itemize} A \emph{viscosity solution} is a continuous
function which is both a subsolution and a supersolution.
\end{de}
One of the motivation for these relaxed boundary conditions is the
stability under uniform convergence. Actually, if the domain $\Om$
satisfies the exterior sphere condition and $F$ is uniformly
elliptic, viscosity subsolutions (resp., supersolutions) satisfy
in the viscosity sense $\frac{\p u}{\p \overrightarrow{n}}\leq $
(resp., $\geq$ )$g(x,u)$  for any $x\in\partial\Om$, see e.g.
Proposition 2.1 in \cite{p2}.

We assume throughout the paper that $\Om$ is a bounded domain of
$\Rn$ of class $C^2$.
\begin{thm}[Strong Comparison Principle, \cite{p2} Theorem 3.1]\label{stcompneu} Assume
that $c$ and $f$ are continuous on $\Oms$. Let $u\in USC(\Oms)$
and $v\in LSC(\Oms)$ be respectively a sub and a supersolution of
\begin{equation*}
\begin{cases}
 \M(D^2u)+c(x)u= f(x) & \text{in} \quad\Om, \\
 \frac{\p u}{\p
\overrightarrow{n}} =\al u & \text{on} \quad\partial\Om. \\
 \end{cases}
 \end{equation*}
If $u\leq v$ on $\Oms$ then either $u<v$ on $\Oms$ or $u\equiv v$
on $\Oms$.
\end{thm}

\begin{prop}[Maximum Principle for $\lam<\lams$, \cite{p2} Theorem 4.5]\label{maxprinc} Assume
$\lam<\lams$. Let $v\in USC(\Oms)$ be a viscosity subsolution of
\begin{equation}\label{sysmaxp}
\begin{cases}
 \M(D^2v)+\lam v=0 & \text{in} \quad\Om, \\
 \frac{\p v}{\p
\overrightarrow{n}} =\al v & \text{on} \quad\partial\Om, \\
 \end{cases}
 \end{equation}
then $v\leq 0$ on $\Oms$.
\end{prop}
\begin{prop}[Minimum Principle for $\lam<\lamso$, \cite{p2} Remark 4.6]\label{minprinc} Assume
$\lam<\lamso$. Let $v\in LSC(\Oms)$ be a viscosity supersolution
of
\begin{equation*}
\begin{cases}
 \M(D^2v)+\lam v=0 & \text{in} \quad\Om, \\
 \frac{\p v}{\p
\overrightarrow{n}} =\al v & \text{on} \quad\partial\Om, \\
 \end{cases}
 \end{equation*}
then $v\geq 0$ on $\Oms$.
\end{prop}

\section{Lipschitz estimates}
In this section we shall prove a local Lipschitz regularity result
for solutions of the Neumann problem associated to general
uniformly elliptic operators, that we will use in the next
sections.
Let us consider the Neumann problem
\begin{equation}\label{lipgensys}
\begin{cases}
 F(x,u,Du,D^2u)=  f(x) & \text{in} \quad\Om, \\
 \frac{\p u}{\p\overrightarrow{n}}= g(x) & \text{on} \quad\partial\Om, \\
 \end{cases}
 \end{equation}
where the operator $F$ is supposed to be continuous on
$\overline{\Om}\times\R\times\R^n \times \emph{S(n)}$ and
satisfying the following assumptions:
\renewcommand{\labelenumi}{(F\arabic{enumi})}
\begin{enumerate}
 \item There exist $b,c>0$ such that for $x\in\overline{\Om},\,r,s\in\R,\, p,q\in \R^n,\, X, Y\in \emph{S(n)}$
 \begin{equation*}\begin{split} \mathcal{M}_{a,A}^-(Y-X)-b|p-q|-c|r-s|&\leq
 F(x,r,p,Y)-F(x,s,q,X)\\&
 \leq \mathcal{M}_{a,A}^+(Y-X)+ b|p-q|+c|r-s|.\end{split}\end{equation*}
 \item There exists $C_1>0$ such that for all $x,y\in\Oms$ and $X\in\emph{S(n)}$
 $$|F(x,0,0,X)-F(y,0,0,X)|\leq C_1|x-y|^\frac{1}{2}\|X\|.$$
\end{enumerate}
\begin{prop}\label{regolaritaloc} Assume that (F1) and (F2)
 hold. Let   $f:\Oms\rightarrow\R$  be  bounded,
$g:\p\Om\rightarrow\R$ be Lipschitz continuous. Let $u\in C(\Oms)$
be a viscosity solution of \eqref{lipgensys}, then, for any
$x_0\in\Oms$ and for any $\rho>0$, there exists $K>0$ such that
 \beq\label{uliploc}|u(x)-u(y)|\leq (MK+|g|_{L^\infty(\p\Om)})|x-y|\quad \forall
x,y\in B_\rho(x_0)\cap\Oms\eeq and
 \beq\label{klipestlocal}\begin{split} K^2-b K&
 \leq C\left[c|u|_{L^\infty(\overline{B}_{3\rho}(x_0)\cap\Oms)}
 +|f|_{L^\infty(\Oms)}
 +\left(1+b\right)|g|_{C^{0,1}(\p\Om)}+\frac{b}{\rho}+\frac{1}{\rho^2}+1\right],\end{split}\eeq
where $M\leq
C(|u|_{L^\infty(\overline{B}_{3\rho}(x_0)\cap\Oms)}+|g|_{L^\infty(\p\Om)}+1)$
  and $C$ depends on $a,\,A,\,C_1,\,n$ and $\Om$.
 \end{prop}


\begin{cor}\label{corregu}Assume $\lam\in\R$ and $\al\geq0$. Let $u\in
C(\Oms)$ be a viscosity solution of
\begin{equation}\label{generallameq}
\begin{cases}
 \M(D^2u)+\lam u=f(x) & \text{in} \quad\Om, \\
 \frac{\p u}{\p
\overrightarrow{n}} =\al u & \text{on} \quad\partial\Om. \\
  \end{cases}
 \end{equation}
 Then, for any $\rho>0$, there exists $K>0$ such that for any $x,y\in\Om_\rho:=\{x\in\Oms\,|\,d(x)\leq\rho\}$
 $$|u(x)-u(y)|\leq (\al|e^{\al
 d(x)}u|_{L^\infty(\Om_\rho)}+MK)|x-y|$$ and
 \beq\label{kvlocal}K^2-C\al K\leq C\left[\left(\al+\al^2+|\lam|\right)|e^{\al
 d(x)}u|_{L^\infty(\Om_{3\rho})}+ |e^{\al
 d(x)}f|_{L^\infty(\Oms)}+\frac{\al}{\rho}+\frac{1}{\rho^2}+1\right],\eeq where $M\leq C(|e^{\al
 d(x)}u|_{L^\infty(\Om_{3\rho})}+1)$  and $C$ depends on $a,A,n$ and $\Om$.
 \end{cor}

 \dims{\bf of Proposition \ref{regolaritaloc}}
 We follow the proof of Proposition III.1 of \cite{il}, that we modify taking test functions which depend on the distance function and that are suitable
 for the Neumann boundary conditions. Moreover, as in \cite{bdl}, we are going to use a regularization of
 $g$. In order to do so, it is convenient to introduce the
 following classical lemma.
 \begin{lem}\label{gextension}Assume $\rho\in C^\infty(\R^n)$, $\rho>0$,
 supp$(\rho)\subset B_1(0)$ and $\int_{\R^n}\rho(y)dy=1$. If $g\in
 C^{0,1}(\R^n)$ and $g$ is bounded, then the function
 $\widetilde{g}:\R^n\times[0,+\infty)\rightarrow\R$ defined by
 \beqs\widetilde{g}(x,\ep):=\int_{\R^n}g(z)\rho\left(\frac{x-z}{\ep}\right)\frac{1}{\ep^n}dz,\,\ep>0,\eeqs
 \beqs \widetilde{g}(x,0)=g(x)\quad\text{for }x\in\R^n,\eeqs is in
 $C^{0,1}(\R^n\times[0,+\infty))$. Moreover, the function
 $\widetilde{g}$ is in $C^2(\R^n\times(0,+\infty))$ with
 \beqs|D_x\widetilde{g}(x,\ep)|,\,|D_\ep\widetilde{g}(x,\ep)|\leq
 C_0,\eeqs
 \beqs
 |D^2_{xx}\widetilde{g}(x,\ep)|,\,|D^2_{x\ep}\widetilde{g}(x,\ep)|,\,|D^2_{\ep\ep}\widetilde{g}(x,\ep)|\leq
 \frac{C_0}{\ep}\quad\text{in }\R^n\times(0,+\infty)\eeqs for some
 positive constant $C_0\leq C|g|_{C^{0,1}(\Rn)}$, with $C$  depending only on $\rho$ and $n$.
 \end{lem}

We first extend $g$ to a $C^{0,1}$ function of $\R^n$ and we still
denote by $g$ this extension. Then, we consider the function
$\widetilde{g}$ associated to $g$ as in Lemma \ref{gextension}.

Since $\Om$ is a domain of class $C^2$, it satisfies the uniform
exterior sphere condition, i.e., there exists $r>0$ such that
$B(x+r\overrightarrow{n}(x),r)\cap \Om =\emptyset$  for any $x\in
\partial\Om.$ From this property it follows that
\begin{equation}\label{sferaest}
\langle
\overrightarrow{n}(x),y-x\rangle\leq\frac{1}{2r}|y-x|^2\quad\text{for
}x\in\partial\Om \text{ and } y\in\Oms.\end{equation} Moreover,
the $C^2$-regularity of $\Om$ implies the existence of a
neighborhood of $\partial\Om$ in $\Oms$ on which the distance from
the boundary
 $$d(x):=\inf\{|x-y|, y\in\partial\Om\},\quad x\in\Oms$$ is of class $C^2$. We still denote by $d$ a
$C^2$ extension of the distance function to the whole $\Oms$.
Without loss of generality we can assume that $|Dd(x)|\leq1$ on
$\Oms$.

Let  $x_0\in\Oms$ and $\rho>0$. Let us denote
$\overline{B}_{\Oms}(x_0,\rho):=\overline{B}_\rho(x_0)\cap\Oms$
and $B_{\Oms}(x_0,\rho):=B_\rho(x_0)\cap\Oms$.
 We are going to show that $u$ is Lipschitz continuous on
$\overline{B}_{\Oms}(x_0,\rho)$. For this aim, let us introduce
the functions
$$\Phi(x)=MK|x|-M(K|x|)^2,$$
$$\Psi_1(x,y)=e^{-L(d(x)+d(y))}\Phi(x-y),$$
$$\Psi_2(x,y)=\widetilde{g}\left(\frac{x+y}{2},
(\delta^2+|x-y|^2)^\frac{1}{2}\right)(d(x)-d(y)),$$ and
$$\varphi (x,y)=\Psi_1(x,y)-\Psi_2(x,y),$$ where
$L$ is a fixed number greater than $\frac{1}{r}$ with $r$ the
radius in \eqref{sferaest}, $K$ and $M$ are  positive constants to
be chosen later and $\delta$ is a small parameter. We also use the
notation
$$\gt(Z,T)=\widetilde{g}\left(\frac{x+y}{2},
(\delta^2+|x-y|^2)^\frac{1}{2}\right).$$ If $K|x|\leq
\frac{1}{4}$, then
\begin{equation}\label{phimagg}\Phi(x)\geq
\frac{3}{4}MK|x|.\end{equation} We define
$$\Delta _K:=\left\{(x,y)\in \R^n\times\R^n:\, |x-y|\leq\frac{1}{4K}\right\}.$$ We
fix $M>1$ and $j>0$ such that
\begin{equation}\label{M}\max_{\overline{B}_{\Oms}(x_0,\rho)^{\,2}}|u(x)-u(y)|+2d_0(|g|_\infty+C_0\delta)\leq
e^{-2Ld_0}\frac{M}{8},\end{equation} $$j=\frac{M}{\rho^2},$$ where
$d_0=\max_{x\in\Oms}d(x),$ and we claim that taking $K$ large
enough, for any small $\delta$ one has
\begin{equation}\label{u-vindelta}
 u(x)-u(y)-\varphi(x,y)-je^{-Ld(x)}|x-x_0|^2\leq 0\quad\text{for }(x,y)\in
\Delta_K\cap\overline{B}_{\Oms}(x_0,\rho)^2.
\end{equation}
 To show \eqref{u-vindelta} we suppose by contradiction that the
maximum of $u(x)-u(y)-\varphi(x,y)-je^{-Ld(x)}|x-x_0|^2$ on
$\Delta_K\cap\overline{B}_{\Oms}(x_0,\rho)^2$ is positive. Then,
for $\delta$ small enough, there is $(\xs,\ys)\in\Delta_K\cap
\overline{B}_{\Oms}(x_0,\rho)^2$ such that $\xs\neq \ys$ and
\begin{equation}\label{u-vcontr}\begin{split} u(\xs)- u(\ys)-\widetilde{\varphi}(\xs,\ys)
=\max_{\Delta_K\cap \overline{B}_{\Oms}(x_0,\rho)\,^2}( u(x)-
u(y)-\widetilde{\varphi}(x,y))>0,\end{split}\end{equation} where
$$\widetilde{\varphi}(x,y)=\varphi(x,y)+je^{-Ld(x)}|x-x_0|^2-C_0\delta(d(x)+d(y)),$$with
$C_0$ the constant defined as in Lemma \ref{gextension}.

The point $(\xs,\ys)$ belongs to $\text{int}(\Delta_K)\cap
B_{\Oms}(x_0,\rho)^2$. Indeed, if $|x-y|=\frac{1}{4K}$, by
\eqref{M} and \eqref{phimagg}, we have
\begin{equation*}\begin{split}
 u(x)- u(y)&\leq  e^{-2Ld_0}
\frac{M}{8}-2d_0(|g|_\infty+C_0\delta)\leq
e^{-L(d(x)+d(y))}\frac{1}{2}MK|x-y|-\Psi_2(x,y)\\&-C_0\delta(d(x)+d(y))\leq\widetilde{\varphi}(x,y).\end{split}\end{equation*}
On the other hand, if $|x-x_0|=\rho$, then
\begin{equation*}\begin{split}
 u(x)- u(y)&\leq e^{-Ld_0}M-2d_0(|g|_\infty+C_0\delta)\leq
 e^{-Ld(x)}\frac{M}{\rho^2}|x-x_0|^2+\Psi_1(x,y)-\Psi_2(x,y)\\&-C_0\delta(d(x)+d(y))=\widetilde{\varphi}(x,y).\end{split}\end{equation*}
Similarly, if $|y-x_0|=\rho$ and $K>K_0/\rho$, for some constant
$K_0>0$, then $ u(x)- u(y)\leq \widetilde{\varphi}(x,y)$.
 Hence, $(\xs,\ys)\in\text{int}(\Delta_K)\cap
B_{\Oms}(x_0,\rho)^2$.

Since $\xs\neq \ys$ we can compute the derivatives of
$\widetilde{\varphi}$ at $(\xs,\ys)$ obtaining
\begin{equation*}\begin{split}D_x\widetilde{\varphi}(\xs,\ys)=&-Le^{-L(d(\xs)+d(\ys))}MK|\xs-\ys|(1-K|\xs-\ys|)Dd(\xs)\\&+
e^{-L(d(\xs)+d(\ys))}MK(1-2K|\xs-\ys|)\frac{(\xs-\ys)}{|\xs-\ys|}-C_0\delta
Dd(\xs)\\&-jLe^{-Ld(\xs)}|\xs-x_0|^2Dd(\xs)+2je^{-Ld(\xs)}(\xs-x_0)-D_x\Psi_2(\xs,\ys),\end{split}\end{equation*}
\begin{equation*}\begin{split}D_y\widetilde{\varphi}(\xs,\ys)=&-Le^{-L(d(\xs)+d(\ys))}MK|\xs-\ys|(1-K|\xs-\ys|)Dd(\ys)\\&-e^{-L(d(\xs)+d(\ys))}MK
(1-2K|\xs-\ys|)\frac{(\xs-\ys)}{|\xs-\ys|}-C_0\delta
Dd(\ys)\\&-D_y\Psi_2(\xs,\ys),\end{split}\end{equation*}where
\begin{equation*}\begin{split}D_x\Psi_2(\xs,\ys)&=\frac{d(\xs)-d(\ys)}{2}D_Z\widetilde{g}\left(\frac{\xs+\ys}{2},
(\delta^2+|\xs-\ys|^2)^\frac{1}{2}\right)\\&+(d(\xs)-d(\ys))\frac{(\xs-\ys)}{(\delta^2+|\xs-\ys|^2)^\frac{1}{2}}D_T\widetilde{g}\left(\frac{\xs+\ys}{2},
(\delta^2+|\xs-\ys|^2)^\frac{1}{2}\right)
\\&+\widetilde{g}\left(\frac{\xs+\ys}{2},
(\delta^2+|\xs-\ys|^2)^\frac{1}{2}\right)Dd(\xs)\end{split}\end{equation*}and
\begin{equation*}\begin{split}D_y\Psi_2(\xs,\ys)&=\frac{d(\xs)-d(\ys)}{2}D_Z\widetilde{g}\left(\frac{\xs+\ys}{2},
(\delta^2+|\xs-\ys|^2)^\frac{1}{2}\right)\\&-(d(\xs)-d(\ys))\frac{(\xs-\ys)}{(\delta^2+|\xs-\ys|^2)^\frac{1}{2}}D_T\widetilde{g}\left(\frac{\xs+\ys}{2},
(\delta^2+|\xs-\ys|^2)^\frac{1}{2}\right)
\\&-\widetilde{g}\left(\frac{\xs+\ys}{2},
(\delta^2+|\xs-\ys|^2)^\frac{1}{2}\right)Dd(\ys).\end{split}\end{equation*}

Observe that
\begin{equation}\label{dxphistima}
|D_x\widetilde{\varphi}(\xs,\ys)|,|D_y\widetilde{\varphi}(\xs,\ys)|\leq
C(MK+C_0+j\rho).\end{equation} Here and henceforth C denotes
various positive constants independent of $K,b,c,f,g$ and $u$.

By Lemma \ref{gextension} \beqs
\left|\widetilde{g}\left(\frac{\xs+\ys}{2},
(\delta^2+|\xs-\ys|^2)^\frac{1}{2}\right)-g(\xs)\right|\leq
C_0(2|\xs-\ys|+\delta),\eeqs then, if $\xs\in\partial \Om$ we have
\begin{equation*}\begin{split}-\langle\overrightarrow{n}(\xs),D_x\Psi_2(\xs,\ys)\rangle-g(\xs)\geq
-C_0(4|\xs-\ys|+\delta).\end{split}\end{equation*}

Hence, using \eqref{sferaest}, if $\xs\in\partial \Om$ we get
\begin{equation}\label{condbordo}\begin{split}
\langle
\overrightarrow{n}(\xs),D_x\widetilde{\varphi}(\xs,\ys)\rangle
-g(\xs)&=Le^{-Ld(\ys)}MK|\xs-\ys|(1-K|\xs-\ys|)
\\&+e^{-Ld(\ys)}MK(1-2K|\xs-\ys|)\langle
\overrightarrow{n}(\xs),\frac{(\xs-\ys)}{|\xs-\ys|}\rangle \\&
+jL|\xs-x_0|^2+2j\langle \overrightarrow{n}(\xs),\xs-x_0\rangle\\&
-\langle\overrightarrow{n}(\xs),D_x\Psi_2(\xs,\ys)\rangle-g(\xs)+C_0\delta\\&\geq
\frac{1}{2}e^{-Ld(\ys)}MK|\xs-\ys|\left(\frac{3}{2}L-\frac{1}{r}\right)+j|\xs-x_0|^2\left(L-\frac{1}{r}\right)\\&-4C_0|\xs-\ys|>0,
\end{split}\end{equation}
for $MK>\frac{16rC_0e^{Ld_0}}{(3rL-2)}$, since $\xs\neq\ys$ and
$L>\frac{1}{r}$. Similarly, if $\ys\in\partial\Om$ then
\begin{equation*}
\langle
\overrightarrow{n}(\ys),-D_y\widetilde{\varphi}(\xs,\ys)\rangle
-g(\ys)\leq
\frac{1}{2}e^{-Ld(\xs)}MK|\xs-\ys|\left(-\frac{3}{2}L+\frac{1}{r}\right)+4C_0|\xs-\ys|<0.\end{equation*}
Then, by definition of sub and supersolution
$$F(\xs, u(\xs),D_x\widetilde{\varphi}(\xs,\ys),X)\geq  f(\xs),\quad\text{if }(D_x\widetilde{\varphi}(\xs,\ys),X)\in \overline{J}^{2,+}
u(\xs),$$ $$F(\ys, u(\ys),-D_y\widetilde{\varphi}(\xs,\ys),Y)\leq
 f(\ys)\quad\text{if }(-D_y\widetilde{\varphi}(\xs,\ys),Y)\in
\overline{J}^{2,-} u(\ys).$$ Since
$(\xs,\ys)\in\text{int}\Delta_K\cap B_{\Oms}(x_0,\rho)\,^2$, it is
a local maximum point of $ u(x)- u(y)-\widetilde{\varphi}(x,y)$ in
$\Oms\,^2$. Then applying Theorem 3.2 in \cite{cil}, for every
$\epsilon>0$ there exist $X,Y\in \emph{S(n)}$ such that $$
(D_x\widetilde{\varphi}(\xs,\ys),X-C_0\delta
D^2d(\xs)+D^2(je^{-Ld(x)}|x-x_0|^2))\in \overline{J}\,^{2,+}
u(\xs),$$ $$(-D_y\widetilde{\varphi}(\xs,\ys),Y+C_0\delta
D^2d(\ys))\in \overline{J}\,^{2,-} u(\ys)$$ and
\begin{equation}\label{tm2ishii}\begin{split}\left(%
\begin{array}{cc}
  X & 0 \\
  0 & -Y \\
\end{array}%
\right)&\leq D^2\varphi(\xs,\ys)+\epsilon
(D^2\varphi(\xs,\ys))^2\\&\leq
D^2\Psi_1(\xs,\ys)-D^2\Psi_2(\xs,\ys)+2\ep(D^2\Psi_1(\xs,\ys))^2+2\ep(D^2\Psi_2(\xs,\ys))^2.
\end{split}\end{equation} Now we want to estimate the matrix on the
right-hand side of the last inequality.

Using Lemma \ref{gextension}, it is easy to check that
\beq\label{psi2hessian}-CC_0\left(\begin{array}{cc}
  I & 0 \\
  0 & I \\
\end{array}%
\right)\leq D^2\Psi_2(\xs,\ys)\leq CC_0\left(\begin{array}{cc}
  I & 0 \\
  0 & I \\
\end{array}%
\right).\eeq Next, let us estimate $D^2\Psi_1(\xs,\ys)$.
\begin{equation*}\begin{split}D^2\Psi_1(\xs,\ys)&=\Phi(\xs-\ys)D^2(e^{-L(d(\xs)+d(\ys))})+D(e^{-L(d(\xs)+d(\ys))})\otimes
D(\Phi(\xs-\ys))\\&+D(\Phi(\xs-\ys))\otimes
D(e^{-L(d(\xs)+d(\ys))})+e^{-L(d(\xs)+d(\ys))}D^2(\Phi(\xs-\ys)).\end{split}\end{equation*}We
set $$A_1:=\Phi(\xs-\ys)D^2(e^{-L(d(\xs)+d(\ys))}),$$
$$A_2:=D(e^{-L(d(\xs)+d(\ys))})\otimes
D(\Phi(\xs-\ys))+D(\Phi(\xs-\ys))\otimes
D(e^{-L(d(\xs)+d(\ys))}),$$
$$A_3:=e^{-L(d(\xs)+d(\ys))}D^2(\Phi(\xs-\ys)).$$Observe that
\begin{equation}\label{a1}A_1\leq CMK|\xs-\ys|\left(\begin{array}{cc}
  I & 0 \\
  0 & I \\
\end{array}%
\right).\end{equation}

For $A_2$ we have the following estimate
\begin{equation}\label{a2}A_2\leq
CMK\left(%
\begin{array}{cc}
  I & 0 \\
  0 & I \\
\end{array}%
\right)
+CMK\left(%
\begin{array}{cc}
  I & -I \\
  -I & I \\
\end{array}%
\right)\leq CMK\left(%
\begin{array}{cc}
  I & 0 \\
  0 & I \\
\end{array}%
\right).\end{equation} Indeed for $\xi,\,\eta\in\R^n$ we compute
\begin{equation*}\begin{split}\langle
A_2(\xi,\eta),(\xi,\eta)\rangle&=2Le^{-L(d(\xs)+d(\ys))}\{\langle
Dd(\xs)\otimes D\Phi(\xs-\ys)(\eta-\xi),\xi\rangle\\&+\langle
Dd(\ys)\otimes D\Phi(\xs-\ys)(\eta-\xi),\eta\rangle\}\\&\leq
CMK(|\xi|+|\eta|)|\eta-\xi|\\&\leq
CMK(|\xi|^2+|\eta|^2)+CMK|\eta-\xi|^2.\end{split}\end{equation*}
Now we consider $A_3$. The matrix $D^2(\Phi(\xs-\ys))$ has the
form
$$D^2(\Phi(\xs-\ys))=\left(%
\begin{array}{cc}
  D^2\Phi(\xs-\ys) & - D^2\Phi(\xs-\ys) \\
  - D^2\Phi(\xs-\ys) &  D^2\Phi(\xs-\ys) \\
\end{array}%
\right),$$and the Hessian matrix of $\Phi(x)$ is
\begin{equation}\label{hessianphi}D^2\Phi(x)=\frac{MK}{|x|}\left(I-\frac{x\otimes
x}{|x|^2}\right)-2MK^2I.\end{equation} If we choose
\begin{equation}\label{epsilon}\epsilon=\frac{|\xs-\ys|}{4MKe^{-L(d(\xs)+d(\ys))}},\end{equation}
then we have the following estimates
$$\epsilon A_1^2\leq
CMK|\xs-\ys|^3I_{2n},\quad \epsilon A_2^2\leq
CMK|\xs-\ys|I_{2n},$$
\begin{equation}\label{aprodotti}\begin{split} \epsilon (A_1A_2+A_2A_1)\leq
CMK|\xs-\ys|^2I_{2n},\end{split}\end{equation}
$$\epsilon (A_1A_3+A_3A_1)\leq
CMK|\xs-\ys|I_{2n},\quad \epsilon (A_2A_3+A_3A_2)\leq CMKI_{2n},$$
 where $I_{2n}:=\left(%
\begin{array}{cc}
  I & 0 \\
  0 & I \\
\end{array}%
\right)$. Then using \eqref{psi2hessian}, \eqref{a1}, \eqref{a2},
\eqref{aprodotti} and observing that $$(D^2(\Phi(\xs-\ys)))^2=\left(%
\begin{array}{cc}
  2(D^2\Phi(\xs-\ys))^2 & - 2(D^2\Phi(\xs-\ys))^2 \\
  - 2(D^2\Phi(\xs-\ys))^2 &  2(D^2\Phi(\xs-\ys))^2 \\
\end{array}%
\right),$$from \eqref{tm2ishii} we can conclude that
$$\left(%
\begin{array}{cc}
  X & 0 \\
  0 & -Y \\
\end{array}%
\right)\leq (MO(K)+CC_0)\left(%
\begin{array}{cc}
  I & 0 \\
  0 & I \\
\end{array}%
\right)+\left(%
\begin{array}{cc}
  B & -B \\
  -B & B \\
\end{array}%
\right),$$ where
\begin{equation}\label{matriceB}B=e^{-L(d(\xs)+d(\ys))}\left[D^2\Phi(\xs-\ys)+\frac{|\xs-\ys|}{MK}
(D^2\Phi(\xs-\ys))^2\right].\end{equation}The last inequality can
be rewritten as follows
$$\left(%
\begin{array}{cc}
  \widetilde{X} & 0 \\
  0 & -\widetilde{Y} \\
\end{array}%
\right)\leq\left(%
\begin{array}{cc}
  B & -B \\
  -B & B \\
\end{array}%
\right),$$ with $\widetilde{X}=X-(MO(K)+CC_0)I$ and
$\widetilde{Y}=Y+(MO(K)+CC_0)I.$

Now we want to get a good estimate for
tr($\widetilde{X}-\widetilde{Y}$), as in \cite{il}. For that aim
let
$$0\leq P:=\frac{(\xs-\ys)\otimes (\xs-\ys)}{|\xs-\ys|^2}\leq I.$$
Since $\widetilde{X}-\widetilde{Y}\leq 0$ and
$\widetilde{X}-\widetilde{Y}\leq 4B,$ we have
$$\text{tr}(\widetilde{X}-\widetilde{Y})\leq \text{tr}(P(\widetilde{X}-\widetilde{Y}))\leq 4 \text{tr}(PB).$$ We
have to compute tr($PB$). From \eqref{hessianphi}, observing that
the matrix $(1/|x|^2)x\otimes x$ is idempotent, i.e.,
$[(1/|x|^2)x\otimes x]^2=(1/|x|^2)x\otimes x$, we compute
$$(D^2\Phi(x))^2=\frac{M^2K^2}{|x|^2}(1-4K|x|)\left(I-\frac{x\otimes
x}{|x|^2}\right)+4M^2K^4I.$$ Then, since $\text{tr}P=1$ and
$4K|\xs-\ys|\leq1$, we have
\begin{equation*}\begin{split}\text{tr}(PB)&=e^{-L(d(\xs)+d(\ys))}MK^2(-2+4K|\xs-\ys|)
\leq -e^{-L(d(\xs)+d(\ys))}MK^2<0.
\end{split}\end{equation*}This gives
\beq\label{estimatex-y}|\text{tr}(\widetilde{X}-\widetilde{Y})|=-\text{tr}(\widetilde{X}-\widetilde{Y})\geq
4e^{-L(d(\xs)+d(\ys))}MK^2\geq CMK^2.\eeq
  Since $\|B\|\leq \frac{CMK}{|\xs-\ys|},$ we have
\begin{equation*}\begin{split}\|B\|^{\frac{1}{2}}|\text{tr}(\widetilde{X}-\widetilde{Y})|^{\frac{1}{2}}&\leq
\left(\frac{CMK}{|\xs-\ys|}\right)^{\frac{1}{2}}|\text{tr}(\widetilde{X}-\widetilde{Y})|^{\frac{1}{2}}
\leq
\frac{C}{K^{\frac{1}{2}}|\xs-\ys|^{\frac{1}{2}}}|\text{tr}(\widetilde{X}-\widetilde{Y})|.
\end{split}\end{equation*}The Lemma III.I in \cite{il} ensures the
existence of a universal constant $C$ depending only on $n$ such
that $$\|\widetilde{X}\|, \|\widetilde{Y}\|\leq
C\{|\text{tr}(\widetilde{X}-\widetilde{Y})|+\|B\|^{\frac{1}{2}}|\text{tr}(\widetilde{X}-\widetilde{Y})|^{\frac{1}{2}}\}.$$
Thanks to the above estimates we can conclude that
\begin{equation}\label{normaxy}\|\widetilde{X}\|,\,\|\widetilde{Y}\|\leq
C|\text{tr}(\widetilde{X}-\widetilde{Y})|\left(1+\frac{1}{K^{\frac{1}{2}}|\xs-\ys|^{\frac{1}{2}}}\right).\end{equation}

Now, using assumptions (F1) and (F2) concerning $F$, the
definition of $\widetilde{X}$ and $\widetilde{Y}$ and the fact
that $ u$ is sub and supersolution we compute
\begin{equation*}\begin{split}f(\xs)&\leq
F(\xs,u(\xs),D_x\widetilde{\varphi}(\xs,\ys),X-C_0\delta
D^2d(\xs)+D^2(je^{-Ld(x)}|x-x_0|^2))\\&\leq
F(\xs,u(\xs),D_x\widetilde{\varphi}(\xs,\ys),\widetilde{X})+MO(K)+CC_0+Cj\\&
\leq F(\xs,u(\ys),-D_y\widetilde{\varphi}(\xs,\ys),\widetilde{Y})+
c|u(\xs)-u(\ys)|
+b|D_x\widetilde{\varphi}(\xs,\ys)+D_y\widetilde{\varphi}(\xs,\ys)|\\&+a\text{tr}(\widetilde{X}-\widetilde{Y})+MO(K)+CC_0+Cj\\&
\leq
F(\ys,u(\ys),-D_y\widetilde{\varphi}(\xs,\ys),\widetilde{Y})+C_1|\xs-\ys|^\frac{1}{2}\|\widetilde{Y}\|+2c|u(\ys)|+2b|D_y\widetilde{\varphi}(\xs,\ys)|
\\&+c|u(\xs)-u(\ys)|
+b|D_x\widetilde{\varphi}(\xs,\ys)+D_y\widetilde{\varphi}(\xs,\ys)|+a\text{tr}(\widetilde{X}-\widetilde{Y})+MO(K)+CC_0+Cj\\&
\leq
f(\ys)+C_1|\xs-\ys|^\frac{1}{2}\|\widetilde{Y}\|+2c|u(\ys)|+2b|D_y\widetilde{\varphi}(\xs,\ys)|+c|u(\xs)-u(\ys)|
\\&+b|D_x\widetilde{\varphi}(\xs,\ys)+D_y\widetilde{\varphi}(\xs,\ys)|+a\text{tr}(\widetilde{X}-\widetilde{Y})+MO(K)+CC_0+Cj.
\end{split}\end{equation*}
 From these inequalities, using
\eqref{dxphistima}, \eqref{normaxy} and \eqref{estimatex-y}, for
$K>\overline{K}$, where $\overline{K}$ is a constant depending
only on $a,A,C_1,n$ and $\Om$, we get
\begin{equation}\label{ultimalemm}\begin{split}
&-2|f|_{L^\infty(\Oms)}-4c|u|_{L^\infty(\overline{B}_{\Oms}(x_0,\rho))}-C|g|_{C^{0,1}(\p\Om)}\\&\leq
Cb|D\widetilde{\varphi}|_\infty+MO(K)+Cj
+C|\text{tr}(\widetilde{X}-\widetilde{Y})|(|\xs-\ys|^{\frac{1}{2}}+K^{-\frac{1}{2}})
+a\text{tr}(\widetilde{X}-\widetilde{Y})
\\&\leq CM\left(-K^2+bK+\frac{b}{\rho}+\frac{1}{\rho^2}\right)+Cb|g|_{C^{0,1}(\p\Om)}.
\end{split}\end{equation}Then, since we have chosen $M>1$, for
$K>\overline{K}$ we  obtain \beq\label{stimKneu}K^2-bK\leq
C\left(|f|_{L^\infty(\Oms)}+c|u|_{L^\infty(\overline{B}_{\Oms}(x_0,\rho))}
+(1+b)|g|_{C^{0,1}(\p\Om)}+\frac{b}{\rho}+\frac{1}{\rho^2}\right),\eeq
and this is a contradiction for $K$ large enough.  This implies
that there exists $K$ satisfying \eqref{stimKneu}, such that
\eqref{u-vindelta} holds true. Next, choosing $x=x_0$,
\eqref{u-vindelta} gives
$$u(x_0)-u(y)\leq \varphi(x_0,y)\quad \forall y\in \overline{B}_{\Oms}(x_0,\rho)\cap\Oms.$$
Repeating the proof in $\overline{B}_{\Oms}(x,2\rho)$ for any
$x\in \overline{B}_{\Oms}(x_0,\rho)$, we finally find the $u$
satisfies \eqref{uliploc} and \eqref{klipestlocal}.\finedim

 \dims {\bf of Corollary \ref{corregu}}
 Let us define
$$v(x):=e^{\al d(x)}u(x).$$ Then, $v$ is a solution of
\begin{equation*}
\begin{cases}
 F(x,v,Dv,D^2v)= e^{\al d(x)}f(x)& \text{in} \quad\Om, \\
 \frac{\p v}{\p
\overrightarrow{n}} =0& \text{on} \quad\partial\Om, \\
  \end{cases}
 \end{equation*} where
 $$F(x,r,p,X)=\M\{X-\al(Dd\otimes p+p\otimes Dd)+\al^2 r(Dd\otimes
 Dd)-\al rD^2d\}+\lam r.$$
 It is easy to check that $F$ satisfies assumptions (F1) and (F2)
 with $C_1=0$, and
 $$c=C(\al^2+\al+|\lam|),\quad b=C\al,$$ where $C$ depends on
 $a,\,A,\,n$ and $\Om$. Then, by Proposition \ref{regolaritaloc}, the
 Lipschitz constant of $v$ on $\Om_\rho$ is bounded from above by $M_vK_v$, where $M_v\leq C(|e^{\al
 d(x)}u|_{L^\infty(\Om_{3\rho})}+1)$  and $K_v$ satisfies
 \eqref{kvlocal}. Hence, for any $x,y\in\Om_\rho$, we have
 \beqs |u(x)-u(y)|\leq|e^{-\al d(x)}-e^{-\al d(y)}||v(x)|+e^{-\al
 d(y)}|v(x)-v(y)|\leq (\al|e^{\al
 d(x)}u|_{L^\infty(\Om_\rho)}+M_vK_v)|x-y|,\eeqs and this concludes
 the proof.\finedim

 \section{Properties of the principal eigenvalues}
 \begin{prop}[Existence of principal eigenfunctions]\label{esistautof}
 There exists $u_\al^+>0$ and $u_\al^-<0$ on $\Oms$ respectively viscosity solution of
\begin{equation}\label{poseigen}
\begin{cases}
 \M(D^2\upl)+\lams \upl=0 & \text{in} \quad\Om, \\
 \frac{\p \upl}{\p
\overrightarrow{n}} =\al \upl & \text{on} \quad\partial\Om, \\
  \end{cases}
 \end{equation}
\begin{equation}\label{negeigen}
\begin{cases}
 \M(D^2\umi)+\lamso\umi=0 & \text{in} \quad\Om, \\
 \frac{\p \umi}{\p
\overrightarrow{n}} =\al \umi & \text{on} \quad\partial\Om. \\
  \end{cases}
 \end{equation}
\end{prop}
\dim We follow the arguments of \cite{bd}. To show the existence
of positive eigenfunctions,  the first step is to prove that if
$f$ is a continuous function such that $f\leq 0$, $f\not\equiv0$,
then for any $\lam<\lams$ there exists a positive solution of
\begin{equation}\label{syslam<lams}
\begin{cases}
 \M(D^2u)+\lam u=f(x) & \text{in} \quad\Om, \\
 \frac{\p u}{\p
\overrightarrow{n}} =\al u & \text{on} \quad\partial\Om. \\
  \end{cases}
 \end{equation}
Observe that $v\equiv1$ is a positive subsolution of
\eqref{sysmaxp} for $\lam\geq0$. This implies, by Proposition
\ref{maxprinc}, that if $\lam<\lams$ then $\lam<0$. Let
$(v_{n})_n$ be the sequence defined by $v_1=0$
 and $v_{n+1}$ be the solution of
\begin{equation*}
\begin{cases}
F(x,v_{n+1},Dv_{n+1},D^2v_{n+1})-(c-\lam) v_{n+1}=e^{\al d(x)}f(x)-cv_n & \text{in} \quad\Om, \\
 \frac{\p v_{n+1}}{\p
\overrightarrow{n}} =0 & \text{on} \quad\partial\Om, \\
  \end{cases}
 \end{equation*}
where $$F(x,r,p,X)=\M\{X-\al(Dd\otimes p+p\otimes Dd)+\al^2
r(Dd\otimes
 Dd)-\al rD^2d\}$$ and $c=C(\al^2+\al)$.  By  comparison, the sequence is positive and increasing. Let $(u_n)_n$ be the sequence defined by
 $u_n(x):=e^{-\al d(x)}v_n(x)$, then $u_{n+1}$ is solution of
\begin{equation*}
\begin{cases}
 \M(D^2u_{n+1})-(c-\lam)u_{n+1}=f(x)-cu_n & \text{in} \quad\Om, \\
 \frac{\p u_{n+1}}{\p
\overrightarrow{n}} =\al u_{n+1} & \text{on} \quad\partial\Om. \\
  \end{cases}
 \end{equation*}
 We claim
that $(u_n)_n$ is  bounded. Suppose that it is not, then defining
$w_n:=\frac{u_n}{|u_n|_\infty}$ one gets that $w_{n+1}$ is a
solution of
\begin{equation*}
\begin{cases}
 \M(D^2w_{n+1})-(c-\lam)w_{n+1}=\frac{f(x)}{|u_{n+1}|_\infty}-c\frac{|u_n|_\infty}{|u_{n+1}|_\infty}w_n & \text{in} \quad\Om, \\
 \frac{\p w_{n+1}}{\p
\overrightarrow{n}} =\al w_{n+1} & \text{on} \quad\partial\Om. \\
  \end{cases}
 \end{equation*}
By Corollary \ref{corregu}, $(w_n)_n$ converges along a
subsequence to a positive function $w$ which
 satisfies
\begin{equation*}
\begin{cases}
\M(D^2w)+\lam w=c(1-k)w\geq0 & \text{in} \quad\Om, \\
 \frac{\p w}{\p
\overrightarrow{n}}=\al w & \text{on} \quad\partial\Om, \\
  \end{cases}
 \end{equation*}
where
 $k:=\limsup_{n\rightarrow+\infty}\frac{|u_n|_\infty}{|u_{n+1}|_\infty}\leq
 1$. This contradicts the Maximum Principle, Proposition
 \ref{maxprinc}. Then $(u_n)_n$ is bounded and letting $n$ go to infinity, by the compactness result, the sequence
converges uniformly to a function $u$ which is a solution of
\eqref{syslam<lams}. Moreover, $u$ is positive by the Strong
Comparison Principle, Theorem \ref{stcompneu}.

We are now in position to construct a sequence $(u_n)_n$ of
positive solutions of
\begin{equation*}
\begin{cases}
\M(D^2u_{n})+\lam_n u_{n}=-1 & \text{in} \quad\Om, \\
 \frac{\p u_{n}}{\p
\overrightarrow{n}} =\al u_n & \text{on} \quad\partial\Om, \\
  \end{cases}
 \end{equation*}
where $(\lam_n)_n$ is an increasing sequence which converges to
$\lams$. The sequence $(u_n)_n$ is unbounded, otherwise one would
 contradict the definition of $\lams$ (see Theorem 8 of
 \cite{bd}).
 Then, up to subsequence, $|u_n|_\infty\rightarrow+\infty$ as
 $n\rightarrow+\infty$ and defining $\phi_n:=\frac{u_n}{|u_n|_\infty}$
 one gets that $\phi_n$ satisfies
\begin{equation*}
\begin{cases}
\M(D^2\phi_{n})+\lam_n \phi_{n}=-\frac{1}{|u_n|_\infty} & \text{in} \quad\Om, \\
 \frac{\p \phi_{n}}{\p
\overrightarrow{n}} =\al \phi_n & \text{on} \quad\partial\Om. \\
  \end{cases}
 \end{equation*}
By Corollary \ref{corregu}, an extracted subsequence converges to
a function $\upl$ with $|\upl|_\infty=1$, which is a solution of
 \eqref{poseigen}. Moreover, by Theorem \ref{stcompneu}, $\upl>0$ on
$\Oms$.

Similar arguments show the existence of negative solutions of
\eqref{negeigen}.

 \finedim

\begin{prop}[Simplicity of the first eigenvalues, \cite{p2} Proposition 7.1]\label{simplicityprop} Let $v\in C(\Oms)$
be a viscosity subsolution (resp. supersolution) of
\eqref{poseigen} (resp. \eqref{negeigen}), then there exists
$t\in\R$ such that $v\equiv t \upl$ (resp. $v\equiv t \umi$).
\end{prop}

\begin{rem}\label{boundprop}{\em Remark that
\beq\label{lamsbound} \lams< -A\al^2,\eeq \beq\label{lamsobound}
\lamso< -a\al^2.\eeq Indeed, the function $v(x):=e^{\al x_1}$,
where $x_1$  is the first coordinate of $x\in\R^n$, is a positive
subsolution of
\begin{equation}\label{propboundeq}
\begin{cases}
 \M(D^2v)-A\al^2 v=0 & \text{in} \quad\Om,\\
 \frac{\p v}{\p
\overrightarrow{n}} =\al v & \text{on} \quad\partial\Om. \\
 \end{cases}
 \end{equation}
Then the Maximum Principle, Proposition \ref{maxprinc}, implies
that $\lams\leq -A\al^2$. If $\lams=-A\al^2$, then by Proposition
\ref{simplicityprop}, $v(x)$ is a solution of \eqref{poseigen} and
this implies that $\Om=\Rn$. Hence \eqref{lamsbound} holds true.
Similarly, inequality \eqref{lamsobound} is a consequence of the
Minimum Principle, Proposition \ref{minprinc}, of Proposition
\ref{simplicityprop} and the fact that $-v(x)$ is a negative
supersolution of \eqref{propboundeq} with $A$ replaced by
$a$.}\end{rem}
\begin{rem}{\em Since $\lams,\lamso<0$ the operator $\M(D^2u)+\lam u$, with $\lam=\lams$
or $\lam=\lamso$ satisfies the Dirichlet Comparison Principle.}\end{rem}
\begin{prop}\label{monot}The sequences $(\lams)_\al$ and $(\lamso)_\al$ are
decreasing.
\end{prop}
\dim Let us prove that $(\lams)_\al$ is decreasing. Consider
$0<\al_1<\al_2$ and let $u_{\al_1}^+$ be a solution of
\eqref{poseigen} with $\al=\al_1$. Then $u_{\al_1}^+$ is a
positive subsolution of
\begin{equation*}
\begin{cases}
 \M(D^2u)+\lam_{\al_1}^+u=0 & \text{in} \quad\Om, \\
 \frac{\p u}{\p
\overrightarrow{n}} =\al_2 u & \text{on} \quad\partial\Om, \\
  \end{cases}
 \end{equation*}
and the Maximum Principle, Proposition \ref{maxprinc}, implies
$\lam_{\al_1}^+\geq \lam_{\al_2}^+$. The strict inequality
$\lam_{\al_1}^+> \lam_{\al_2}^+$ follows from Proposition
\ref{simplicityprop}. \finedim

\begin{lem}\label{maxboundlem} Let $\upl$ and $\umi$ be respectively a positive solution of \eqref{poseigen} and a negative  solution
of \eqref{negeigen}, then \beqs \upl(x)<\max_{\p\Om}
\upl\quad\forall x\in\Om,\eeqs \beqs\umi(x)>\min_{\p\Om}
\umi\quad\forall x\in\Om.\eeqs
\end{lem}
\dim Let us show the result for $\upl$. Suppose by contradiction
that the maximum of $\upl$ is attained at some point $x_0\in\Om$
and let $v(x):=\upl(x)-\upl(x_0)$. Since $\upl(x_0)>0$ and
$\lams<0$, $v$ satisfies
$$\M(D^2 v)+\lams v\geq 0\quad\text{in }\Om$$ and $v\leq 0$ in
$\Om$, $v(x_0)=0$. Then the Strong Maximum Principle implies
$\upl\equiv \upl(x_0)$ in $\Om$ and this contradicts the fact that
$\upl$ solves \eqref{poseigen}. \finedim

\section{Liouville type results}
For $\gamma>0$ let us introduce the system
\beq\label{pbhalfspace}\left\{%
\begin{array}{ll}
    \M(D^2 u)- \gamma u=0 & \hbox{in }\Rp, \\
    -\frac{\p u}{\p x_n}=u & \hbox{on } \p\Rn.\\
\end{array}%
\right.\eeq

\begin{thm}\label{halfspacethm}If $\gamma >A$,  any bounded subsolution of \eqref{pbhalfspace} is non-positive in $\Rp$.

If $\gamma >a$, any bounded supersolution of \eqref{pbhalfspace}
is non-negative in $\Rp$.

Hence, if $\gamma>A$ there are no, non trivial bounded solutions
of \eqref{pbhalfspace}.
\end{thm}
\begin{rem}{\em It turns out that Theorem \ref{halfspacethm} is
sharp:  $u(x)=e^{-x_n}$ (resp., $u(x)=-e^{-x_n}$) is a positive
bounded subsolution (resp., negative bounded supersolution) of
\eqref{pbhalfspace} for every $\gamma\leq A$ (resp., $\gamma\leq
a$).}

{\em Theorem \ref{halfspacethm} also fails without the boundedness
condition. Indeed, $u(x)=e^{\nu\cdot x}$ (resp.,
$u(x)=-e^{\nu\cdot x}$), with $\nu=(\nu_1,...,\nu_{n-1},-1)$,
$|\nu|>1$, is an unbounded subsolution (resp., supersolution) of
\eqref{pbhalfspace} for $A<\gamma\leq A|\nu|^2$ (resp.,
$a<\gamma\leq |\nu|^2a$).}
\end{rem}

We assume that  $u(x)$ is a bounded subsolution of
\eqref{pbhalfspace} with $\gamma>0$, which is positive somewhere.
%
We normalize $u$ so that \beq\label{supequal1}\sup_{\Rp} u=1.\eeq
Then $u$ is a viscosity subsolution of
\beq\label{neumhalfspace}\left\{%
\begin{array}{ll}
    \M(D^2 u)- \gamma u=0 & \hbox{in }\Rp, \\
    -\frac{\p u}{\p x_n}=1 & \hbox{on } \p\Rn.\\
\end{array}%
\right.\eeq

\begin{prop}\label{comparisonhalfspace} Assume $\gamma>0$ and $k\in\R$. Let $u\in USC(\R_+^n)$ and $v\in LSC(\R_+^n)$ be
respectively  bounded viscosity sub and supersolution of
\beq\label{neumhalfspacek}\left\{%
\begin{array}{ll}
    \M(D^2 u)- \gamma u=0 & \hbox{in }\Rp,\\
    -\frac{\p u}{\p x_n}=k & \hbox{on } \p\Rn.\\
\end{array}%
\right.\eeq  Then $u\leq v$ in $\Rp$.
\end{prop}
\dim
 Suppose by contradiction that
$\sup_{\Rp}(u-v)=M>0$. Let $\psi$  be a smooth positive function
with bounded derivatives and such that $\psi(x)\rightarrow+\infty$
as $|x|\rightarrow+\infty$. Let $\chi(x)=\chi(x_n)$ be a smooth
function such that $\chi(x_n)=x_n$ for $|x_n|\leq1$ and
$\chi(x_n)\equiv 0$ for $|x_n|>2$. Let
$$\varphi(x,y)=\frac{j}{2}|x-y|^2-k(x_n-y_n)+\beta\psi(x)-\ep(\chi(x)+\chi(y)).$$
Then, for $\beta$ and $\ep$ small enough and $j>0$, the supremum
of the function $u(x)-v(y)-\varphi(x,y)$
 is greater than $\frac{M}{2}$ and it is reached at some point
$(\xs,\ys)\in \overline{\Rp}\times\overline{\Rp}$.

If $\xs\in\p\Rp$ then, for $\ys\in\Rp$,
\beqs\begin{split}-\p_{x_n}\varphi(\xs,\ys)-k=-j(\xs_n-\ys_n)+k
-\beta\p_{x_n}\psi(\xs)+\ep-k=j\ys_n-\beta\p_{x_n}\psi(\xs)+\ep>0\end{split}\eeqs
for $\ep>\beta|D\psi|_\infty$.

If $\ys\in\p\Rp$ then, for $\xs\in\Rp$
\beqs\p_{y_n}\varphi(\xs,\ys)-k=-j(\xs_n-\ys_n)-\ep=-j\xs_n-\ep<
0.\eeqs  Both inequalities contradict the definition of sub and
supersolution, therefore $\xs,\ys\in\Rp$.

Applying Theorem 3.2 of \cite{cil}, there exist
$X,Y\in\emph{S(n)}$ such that $(D_x\varphi(\xs,\ys),X+\beta
D^2\psi(\xs)-\ep D^2\chi(\xs))\in \overline{J}^{2,+}u(\xs)$,
$(-D_y\varphi(\xs,\ys),Y+\ep D^2\chi(\ys))\in
\overline{J}^{2,-}v(\ys)$ and
 \begin{equation*}-3j\left(%
\begin{array}{cc}
  I & 0 \\
  0 & I \\
\end{array}%
\right)\leq \left(%
\begin{array}{cc}
  X & 0 \\
  0 & -Y \\
\end{array}%
\right)\leq 3j\left(%
\begin{array}{cc}
  I & -I \\
  -I & I \\
\end{array}%
\right).
\end{equation*}
Since $u$ and $v$ are respectively sub and supersolution, we have
\beqs\M(X+\beta D^2\psi(\xs)-\ep D^2\chi(\xs))\geq \gamma
u(\xs),\eeqs \beqs \M(Y+\ep D^2\chi(\ys))\leq \gamma v(\ys).\eeqs
Subtracting the two previous inequalities, using the properties of
Pucci's operators and that
$$u(\xs)-v(\ys)>\frac{M}{2}+\frac{j}{2}|\xs_n-\ys_n|^2-k(\xs_n-\ys_n)-\ep(\chi(\xs)+\chi(\ys))\geq
\frac{M}{2}-\frac{k^2}{2j}-C\ep,$$ we finally get
\beqs\begin{split}
\frac{\gamma}{2}\left(M-\frac{k^2}{j}-C\ep\right)&<\gamma(u(\xs)-v(\ys))\leq
\M(X+\beta D^2\psi(\xs)-\ep D^2\chi(\xs))\\&-\M(Y+\ep
D^2\chi(\ys))\\&\leq \M(X-Y)+\beta
\M(D^2\psi(\xs))-\ep\mathcal{M}_{a,A}^-(
D^2\chi(\xs)+D^2\chi(\ys))\\&\leq \beta
\M(D^2\psi(\xs))-\ep\mathcal{M}_{a,A}^-(
D^2\chi(\xs)+D^2\chi(\ys)).\end{split}\eeqs This is a
contradiction for $\beta$ and $\ep$ small enough and $j$ large.
Then $u\leq v$ in $\Rp$. \finedim

\dims {\bf of Theorem \ref{halfspacethm}} The function
$$v(x)=\sqrt{\frac{A}{\gamma}}e^{-\sqrt{\frac{\gamma}{A}}x_n}$$ is
the bounded viscosity solution of \eqref{neumhalfspace}. Then by
 Proposition \ref{comparisonhalfspace}
$$u(x)\leq \sqrt{\frac{A}{\gamma}}e^{-\sqrt{\frac{\gamma}{A}}x_n},
\quad\text{for any }x\in\Rp.$$ It follows from \eqref{supequal1}
that
$$1=\sup_{\Rp}u\leq \sqrt{\frac{A}{\gamma}},$$ i.e. $\gamma\leq
A$.

Similarly, if $u$ is a negative supersolution of
\eqref{pbhalfspace}, normalized so that $\min_{\Rp}u=-1$, then $u$
is a supersolution of
\beqs\left\{%
\begin{array}{ll}
    \M(D^2 u)- \gamma u=0 & \hbox{in }\Rp, \\
    -\frac{\p u}{\p x_n}=-1 & \hbox{on } \p\Rn,\\
\end{array}%
\right.\eeqs and by comparison
$$u(x)\geq
-\sqrt{\frac{a}{\gamma}}e^{-\sqrt{\frac{\gamma}{a}}x_n}.$$ This
implies $\gamma\leq a$ and Theorem \ref{halfspacethm} is
proved.\finedim

\section{Asymptotic behavior and Proof of Theorem \ref{mainthm}}
We start by the following simple result:
\begin{prop} $\displaystyle \lim_{\alpha\rightarrow 0}\lambda^{\pm}_{\alpha}=0$.
\end{prop}
\dim By Proposition \ref{monot}, $\lams$ increases to some value
$\lambda_0\leq 0$. On the other hand, the sequence of normalized
solutions $(\upl)_\alpha$, by the  Lipschitz estimates Corollary
\ref{corregu}, converges to  $u_0$ a positive solution of
\begin{equation*}
\begin{cases}
 \M(D^2u)+\lambda_0 u=0 & \text{in} \quad\Om, \\
 \frac{\p u}{\p
\overrightarrow{n}}=0 & \text{on} \quad\partial\Om, \\
 \end{cases}
\end{equation*}
which satisfies $|u_0|=1$. Recall that $0$ is the principal
eigenvalue for the Neumann problem. If $\lambda_0<0$, the Maximum
Principle below the first eigenvalue, i.e. Proposition
\ref{maxprinc}, implies that $u_0\leq 0$ a contradiction. \finedim

We consider now the asymptotic behavior at infinity.
By Remark \ref{boundprop}, it is enough to show that
\beq\label{limsuplams}\limsup_{\al\rightarrow+\infty}\frac{\lams}{-\al^2}\leq
A,\eeq and
\beq\label{limsuplamso}\limsup_{\al\rightarrow+\infty}\frac{\lamso}{-\al^2}\leq
a.\eeq

We are going to show \eqref{limsuplams}. For $\al>0$, let $\upl$
be a positive solution of \eqref{poseigen}. By Lemma
\ref{maxboundlem}, we know that $\upl$ attains its maximum at
$x_\al\in\p\Om$. After normalization, we can assume that
$\max_{\Oms}\upl=1$ and $\xa\rightarrow0$ as
$\al\rightarrow+\infty$. Furthermore, we can assume that there is
a $C^2$ function $\phi$ and $r>0$ such that
\begin{eqnarray*}
 &&x_n=\phi(x'),\quad \forall (x',x_n)\in \p\Om\cap B_r(0) \\&&
  x_n>\phi(x'),\quad \forall (x',x_n)\in \Om\cap B_r(0) \\&&
  \phi(0)=0, \\&&
\p_{x_i}\phi(0)=0,\quad\text{for }i=1,...,n-1.\\&&
\end{eqnarray*}
 We flatten  $\p\Om$ near the origin. Let
$\Phi(x):\Om\cap  B_r(0)\rightarrow \Om_{\Phi}:=\Phi(\Om\cap
B_r(0))$, be such that
\beq\label{Phi}\begin{split}&\Phi_i(x)=x_i,\quad i=1,...,n-1,\\&
 \Phi_n(x)=x_n-\phi(x').
\end{split}\eeq

Denote by $x=\Psi(y)$ the inverse of $y=\Phi(x)$. The function
$$v_\al(y)=\upl(\Psi(y))$$ is solution of
\beq\label{systemflatbound}\left\{%
\begin{array}{ll}
    \M\left[ \left(\sum_{l,k=1}^n
    \p^2_{y_ly_k}v_\al\p_{x_j}\Phi_k(\Psi(y))\p_{x_i}\Phi_l(\Psi(y))\right)_{ij}\right.\\
    \quad\left.+\left(\sum_{k=1}^n\p_{y_k}v_\al\p^2_{x_ix_j}\Phi_k(\Psi(y))\right)_{ij}\right]+\lams v_\al=0 & y\in \Om_{\Phi}, \\
    \sum_{k,j=1}^n\p_{y_k}v_\al\p_{x_j}\Phi_k(\Psi(y))\overrightarrow{n}_j(\Psi(y))=\al v_\al  & y\in\p\Om_{\Phi}.\\
\end{array}%
\right.\eeq
 Since the
exterior normal $\overrightarrow{n}(x)$ at $x\in \p\Om\cap B_r(0)$
is
$$\overrightarrow{n}(x)=\frac{(D\phi(x'),-1)}{\sqrt{|D\phi(x')|^2+1}},$$
by \eqref{Phi}, the boundary condition in \eqref{systemflatbound}
can be rewritten as follows
\beqs\frac{1}{\sqrt{|D\phi(y')|^2+1}}\sum_{k=1}^{n-1}\p_{y_k}v_\al\p_{x_k}\phi(y')-\left(\sqrt{|D\phi(y')|^2+1}\right)\p_{y_n}v_\al=\al
v_\al,\quad y\in\p\Om_{\Phi}.\eeqs Notice that, since
$D\phi(x')\rightarrow 0$ as $x'\rightarrow0$,
$D\Phi(\Psi(y))\rightarrow I$ as $y\rightarrow0$, where $I$ is the
identity matrix of $\emph{S(n)}$.

We now consider two different cases.

\noindent\emph{Case 1.}
$$\limsup_{\al\rightarrow+\infty}\frac{\lams}{-\al^2}=\gamma<+\infty.$$
Without loss of generality, we may assume that
$\frac{\lams}{-\al^2}\rightarrow\gamma$ as
$\al\rightarrow+\infty$, and $\upl(\xa)=\max_{\Oms}\upl=1$,
$\xa\rightarrow0$ as $\al\rightarrow+\infty$. We let
$$z=\al(y-y_\al),$$ where $y_\al=\Phi(\xa)=(\xa',0)$. We set
$$w_\al(z)=v_\al(y)=\upl(x),$$ then for any $R>0$, as $\al$ becomes sufficiently large,
$w_\al$ is  solution of
\beq\label{systemblowup}\left\{%
\begin{array}{ll}
    \M\left[ \left(\sum_{l,k=1}^n
    \p^2_{z_lz_k}w_\al\p_{x_j}\Phi_k(\Psi\left(y\right))\p_{x_i}\Phi_l(\Psi\left(y\right))\right)_{ij}\right.\\
    \quad\left.+\frac{1}{\al}\left(\sum_{k=1}^n\p_{z_k}w_\al\p^2_{x_ix_j}\Phi_k(\Psi\left(y\right))\right)_{ij}\right]
    +\frac{\lams}{\al^2}  w_\al=0 & z\in B_{R}^+, \\
    \frac{1}{\sqrt{\left|D\phi\left(y'\right)\right|^2+1}}\sum_{k=1}^{n-1}\p_{z_k}w_\al\p_{x_k}\phi\left(y'\right)
    -\left(\sqrt{\left|D\phi\left(y'\right)\right|^2+1}\right)\p_{z_n}w_\al=w_\al  & z\in \Gamma_{R},\\
\end{array}%
\right.\eeq where $$y=y(z)=\frac{z}{\al}+y_\al$$ and
$$B_R^+:=B_R(0)\cap\Rp,\quad \Gamma_R:=B_R(0)\cap\p\Rp.$$ Since for $z\in
B_{R}^+$, $z/\al+y_\al\rightarrow0$ as $\al\rightarrow+\infty$ and
$\p_{x_i}\phi(0)=0$ for $i=1,...,n-1$, for $\al$ sufficiently
large, $I/2\leq D\Psi\left(z/\al+y_\al\right)\leq 2I$. Hence, if
$L_\al$ is the Lipschitz constant of $\upl$ in the set
$\{x=\Psi(z/\al+y_\al),\,|z|\leq R\}$, we have \beqs
|w_\al(z_1)-w_\al(z_2)|=\left|\upl\left(\Psi\left(\frac{z_1}{\al}+y_\al\right)\right)-\upl\left(\Psi\left(\frac{z_2}{\al}+y_\al\right)\right)\right|
\leq \frac{2L_\al}{\al}|z_1-z_2|.\eeqs

Remark that if $|z|\leq R$, then $d(x)\leq CR/\al$ for
$x=\Psi(z/\al+y_\al)$, where $C$ depends on $\phi$. Hence, since
for  $\rho=CR/\al$, $|e^{\al
d(x)}\upl|_{L^\infty(\Om_{3\rho})}\leq e^{3CR}$,  Corollary
\ref{corregu} gives
$$L_\al\leq C e^{3CR}(\al+K_\al),$$
where $K_\al$ satisfies
$$K_\al^2-C\al K_\al\leq
C\left[(\al+\al^2+|\lams|)e^{3CR}+\frac{\al^2}{CR^2}+1\right].$$
This implies that the sequence $(w_\al)_\al$ is bounded in the
space of Lipschitz continuous functions of $\overline{B}_R^+$ for
any fixed $R>0$, and then, up to subsequence, $w_\al\rightarrow
w_0$ uniformly on $\overline{B}_R^+$, with $\sup_{\Rp}w=1$,
viscosity solution of \eqref{pbhalfspace}. Moreover by the Strong
Comparison Principle, $w>0$ on $\overline{\Rp}$. Then, by Theorem
\ref{halfspacethm}, $\gamma\leq A$ and this proves
\eqref{limsuplams}.
\bigskip

\noindent\emph{Case 2.}
$$\limsup_{\al\rightarrow+\infty}\frac{\lams}{-\al^2}=+\infty.$$

Let $\upl$ be the sequence of positive solutions of
\eqref{poseigen} such that
$$\frac{\lams}{-\al^2}=:l_\al\rightarrow+\infty\quad\text{as
}\al\rightarrow+\infty,$$ and $\upl(\xa)=\max_{\Oms}\upl=1$.
Define
$$z=\sqrt{l_\al}\al(y-y_\al)\quad\text{and}\quad
w_\al(z)=u_\al(x),$$ where $y=\Phi(x)$ and $y_\al=\Phi(x_\al)$.
Then, for any $R>0$, as $\al$ becomes sufficiently large, $w_\al$
satisfies
\beqs\left\{%
\begin{array}{ll}
    \M\left[ \left(\sum_{l,k=1}^n
    \p^2_{z_lz_k}w_\al\p_{x_j}\Phi_k(\Psi(y))\p_{x_i}\Phi_l(\Psi\left(y\right))\right)_{ij}\right.\\
    \quad\left.+\frac{1}{\sqrt{-\lams}}\left(\sum_{k=1}^n\p_{z_k}w_\al\p^2_{x_ix_j}\Phi_k(\Psi\left(y\right))\right)_{ij}\right]
    - w_\al=0 & z\in B_{R}^+, \\
    \frac{1}{\sqrt{\left|D\phi\left(y'\right)\right|^2+1}}\sum_{k=1}^{n-1}\p_{z_k}w_\al\p_{x_k}\phi\left(y'\right)
    -\left(\sqrt{\left|D\phi\left(y'\right)\right|^2+1}\right)\p_{z_n}w_\al=\frac{1}{\sqrt{l_\al}}w_\al  & z\in \Gamma_{R},\\
\end{array}%
\right.\eeqs where $y=y(z)=\frac{z}{\al}+y_\al$. As in Case 1, we
can show that $w_\al\rightarrow w_0$ which is  a bounded positive
viscosity solution of
\beq\label{pureneumann}\left\{%
\begin{array}{ll}
    \M(D^2 w_0)- w_0=0 & \hbox{in }\Rp, \\
    -\frac{\p w_0}{\p x_n}=0 & \hbox{on } \p\Rn.\\
\end{array}%
\right.\eeq On the other hand, by Proposition
\ref{comparisonhalfspace}, the only bounded viscosity solution of
\eqref{pureneumann} is $u\equiv 0$ and we reach a contradiction.
\finedim

\begin{prop}\label{convcompacts}Let $\upl$ and $\umi$ be respectively the normalized
solution of \eqref{poseigen} and \eqref{negeigen}, i.e.
$\|u_\al^{\pm}\|_\infty=1$, then for any compact set $K\subset\Om$
\beqs \|u_\al^{\pm}\|_{L^\infty(K)}\rightarrow0\quad\text{as
}\al\rightarrow+\infty.\eeqs
\end{prop}
\dim Let $\upl$ be the normalized solution of \eqref{poseigen} and
let $K$ be a compact set contained in $\Om$. Let $\xa\in K$ be
such that $\max_K \upl=\upl(\xa)$. Define $z=\al(x-\xa)$ and
$w_\al(z)=\upl(x)$ for $|z|< \al r$ where $r=$dist$(K,\p\Om)$.
Then for any $R>0$, as $\al$ becomes large, $w_\al(z)$ satisfies
\beqs\M(D^2 w_\al)+\frac{\lams}{\al^2} w_\al=0\quad\text{in
}B_{2R},\eeqs and $\|w_\al\|_\infty\leq 1$. By standard elliptic
estimates, see e.g. \cite{cc} and Theorem \ref{mainthm},
$w_\al\rightarrow w_0$ non-negative solution of \beqs \M(D^2 w_0)-
Aw_0=0\quad\text{in }\Rn.\eeqs It is well-know that there are no
nontrivial bounded solutions of the above equation, see e.g.
\cite{cil}, hence $w_\al(0)=\max_K \upl\rightarrow0$ as
$\al\rightarrow+\infty$ and Proposition \ref{convcompacts} is
proved.\finedim

\end{document}